\documentclass [winedt,yap]{iitparc}
\usepackage{cite}
\usepackage{amsmath,amssymb,amsfonts,bm}
\usepackage{lscape}
\usepackage{floatfig,wrapfig,epsfig}
\usepackage{subfigure}
\usepackage{color}
\usepackage{psboxit}
\usepackage{rotating}
\usepackage{curves}
\usepackage{ulem}
\usepackage[mathscr]{eucal}
\RequirePackage{psfrag}
\RequirePackage{graphicx}

\begin{document}\normalem
\initfloatingfigs
\frontmatter          

\IssuePrice{25.00}%
\TransYearOfIssue{2011}%
\TransCopyrightYear{2011}%
\OrigYearOfIssue{2011}%
\OrigCopyrightYear{2011}%

\TransVolumeNo{72}%
\TransIssueNo{3}%
\OrigIssueNo{3}%
\OrigPages{173}


\mainmatter

\setcounter{page}{626}
\CRubrika{NOTES, MEETINGS, BOOK REVIEW}
\Rubrika{NOTES, MEETINGS, BOOK REVIEW}
\def\x#1{} 
\def\cdc{,\ldots,}
\def\la{\lambda}
\def\ind{\mathop{\rm ind}\nolimits}

\OrigJournalName{Avtomatika i Telemekhanika}
\OrigIssueNo{3}
\OrigYearOfIssue{2011}
\OrigCopyrightYear{2011}%

\title{Addendum to the paper ``On determining the eigenprojection and components of a matrix''}%
\titlerunning{Addendum to ``On determining the eigenprojection...''}

\author{P. Yu. Chebotarev and R. P. Agaev}
\authorrunning{Chebotarev \lowercase{and} Agaev}
\OrigCopyrightedAuthors{P.Yu. Chebotarev, R.P. Agaev}

\institute{Trapeznikov Institute of Control Sciences, Russian
Academy of Sciences, Moscow, Russia}
\received{Received November 1, 2010}
\OrigPages{p.~173}

\maketitle

\begin{abstract}
The purpose of this note is to correct an inaccuracy in the paper: 
Agaev,~R.P. and Chebotarev,~P.Yu. On determining the eigenprojection and components of a matrix, {\it Autom.\ Remote Control}, 2002, vol.~63, pp.~1537--1545 and to present one of its results in a simplified form.\\
DOI: 10.1134/S000511791103012X
\end{abstract}


In \cite{AgaChe02AiT}, an error was made which may lead to incorrect signs of the right-hand sides in formulas (16) and~(17) of Proposition~2.
Below we present Proposition~2 in a corrected and simplified form.
\setcounter{proposition}{1}
\setcounter{equation}{15}

\vspace{-2mm}
\begin{proposition}
\label{p-ZZk=}
Let $A\in{\mathbb C}^{n\times n};$ let $Z$ be the eigenprojection of~$A.$ Suppose that $\la_1\cdc\la_s$ are the distinct\x{} eigenvalues of~$A,$ $\nu_1\cdc\nu_s$ are their indices$,$ and the integers $u_1\cdc u_s$ are such that $u_i\ge\nu_i,$\, $i=1\cdc s.$ Let $u\ge\ind A.$ Then$:$
\begin{gather}
\label{e-Zeva}
Z=\prod_{i:\,\la_i\ne0}\Bigl(I-(A/\la_i)^u\Bigr)^{\!u_i}
\end{gather}
\vspace{-1mm} and \vspace{-1mm}
\begin{gather}
\label{e-Zk_eva} Z_{kj}=
\prod_{i\ne k}\left(I-\left(\frac{A-\la_kI}{\la_i-\la_k}\right)^{\!u_k}\right)^{\!u_i}\!(j!)^{-1}(A-\la_kI)^j,
\end{gather}
where $Z_{kj}$ is the order $j$ component of $A$ corresponding to $\la_k,$\, $k=1\cdc s,$\, $j=0\cdc \nu_k-1.$
\end{proposition}
\vspace{-2mm}

If $j=0,$ then Eq.\,\eqref{e-Zk_eva} determines the eigenprojections of $A$ corresponding to its eigenvalues.

\PPR{\ref{p-ZZk=}}{Note that
$\la\,\xi(\la)\equiv\la\prod_{i:\,\la_i\ne0}(\la-\la_i^u)^{u_i}$
is an annihilating polynomial for~$A^u.$ Indeed, it is divisible by the minimal\x{} polynomial for $A^u$ because $\la_i^u$ are the eigenvalues of~$A^u,$ their indices do not exceed the numbers $u_i$ by Lemma~2, and $\ind A^u\le1$ according to~(3). To prove \eqref{e-Zeva}, it suffices now to take $h(\la)$ (see (11)) as the polynomial $\xi(\la)$ divided by its absolute term
$p=\prod_{i:\,\la_i\ne0}(-\la_i^u)^{u_i},\x{}$
apply Theorem~1, , and perform some algebraic transformations.

By applying now \eqref{e-Zeva} to the matrix $A-\la_kI$ whose index, by definition, does not exceed $u_k,$ the distinct eigenvalues $\la_i-\la_k$ have the indices $\nu_i,$ $i=1\cdc s,$ respectively, and so $u_1\cdc u_s$ are the non-strict upper bounds for these indices, we obtain \eqref{e-Zk_eva} with $j=0.$ Now the expressions \eqref{e-Zk_eva} for the components of $A$ of higher orders follow from~(15).\qed}

Proposition~\ref{p-ZZk=} generalizes formula (5.4.3) in~\cite{Lankaster78} which refers to the case of $\nu_1=\ldots=\nu_s=1.$

\end{document}